# A Basic Łukasiewicz *m*-valued conditional logic[①]

Shuquan Huo
School of Philosophy and Public Management,
Henan University
475000 Kaifeng, (P. R. China)
Email: huosq@henu.edu.cn

**Abstract.** This paper is devoted to the construction of conditional logic system of Łukasiewicz *m*-valued propositional logic. We construct conditional logic system **ŁCR** based on Łukasiewicz *m*-valued propositional logic. We construct world semantics for the system by generalizing conditional and accessibility relation from classical bivalent to *m*-valued, and prove its soundness, completeness and finite model property. Conditionals of **ŁCR** cannot be generalized directly to variable strict conditionals, but they are stricter than classical conditionals.

## 1. Introduction

Conditional logics have traditionally been intended to formalize various modes of reasoning and natural language conditionals. The contemporary development of conditional logic can be traced back to the works of Adams [1], Stalnaker [25, 26] and Lewis [16], and a variety of conditional logic systems have been developed with tight connections to modal logic [3, 23] nowadays. Despite the variety, almost all extant conditional logics mainly extend classical logic. As far as we know, only a small number of references in the literature generalize conditional logic to non-classical logic [4, 32]. But we think that it is also possible and necessary to develop the conditional logic of many-valued logic.

The implication of Łukasiewicz logic maintains the properties of classical logic implication, and when the truth value of the antecedent is greater than that of the consequent, its truth value is not necessarily false, but may take a variety of truth values, which is difficult to explain [24]. But we think that since Łukasiewicz only regards the truth value 1 as the designated truth value, truth values other than 1 can be regarded as various forms of falsehood. Thus the implication of Łukasiewicz logic still represents inference or conditionals. However, similar to classical logic, a better representation of conditionals and reasoning requires the development of conditional logic. Many-valued logic can deal with fuzzy phenomena well. Hájek and Novák have proved that fuzzy logic can indeed deal with the sorites paradox [13]. Therefore, many-valued logic has its own

① This research was supported by Major Program of National Fund of Philosophy and Social Science of China (18ZDA032)

special functions, and should have its conditional logic.

Conditional logic is a generalization of modal logic. Many-valued modal logic has been studied in some references in the literature [2, 7, 8, 9, 14, 18, 19, 22, 27, 28, 29]. The modal logic and multi-modal logic of Łukasiewicz many-valued logic are also discussed [6], but we have not found any references in the literature that generalizes conditional logic to many-valued logic yet. Due to the complexity of many-valued accessibility relation, there has not been a Łukasiewicz infinite-valued modal logic whose semantic model has many-valued accessibility relation, so we only consider establishing Łukasiewicz $m$-valued conditional logic. We believe that as long as $m$ is large enough, finite-valued logic is still a kind of fuzzy logic and can still describe fuzzy phenomena.

In [3], Chellas developed the basic (normal) classical conditional logic **CK**, but he put no restrictions on the selection function. In [32], Weiss developed and thoroughly analyzed its intuitionistic counterpart, basic intuitionistic conditional logic **ICK**. Ciardelli and Liu gave more systems of intuitionistic conditional logic [4]. In the model of classical conditional logic, the accessibility relation is bivalent, and the possible world (set) chosen by selection function is the one most similar to the current world (set). According to Jing et al. [15], there can be different degrees of belief, and therefore, accessibility relation can also be many-valued. Lewis argued that counterfactual conditionals are variable strict conditionals [16]. In fact, because the accessibility relation is $m$-valued, it is difficult for us to generalize the constructed conditional directly to the variable strict case. However, Since we take into account not only the case when the accessibility relation is classical, but also other many-valued cases, our conditionals are obviously stricter than the classical conditionals.

The content of this paper is arranged as follows. In Section 2, we will construct conditional logic **ŁCR** on the basis of Łukasiewicz $m$-valued propositional logic, give its Kripke semantics, axiom system and inference rules, and prove its soundness, completeness. We treat a proposition as a $m$-tuple of sets of possible worlds, by which we can solve the difficulty of constructing Łukasiewicz $m$-valued conditional logic in which conditional and accessibility relation are many-valued. In Section 3, we prove the finite model property of **ŁCR**. In Section 4, we construct the extended system of **ŁCR** by slightly expanding. Throughout the paper, we use the technical achievements of Łukasiewicz $m$-valued logic which can define the operators $\mathsf{J}$ and $\mathsf{I}$ to make the proof more convenient [19, 20].

## 2. Basic Łukasiewicz $m$-valued conditional logic ŁCR

### 2.1. Language

The primary aim of this section is to characterize the basic Łukasiewicz $m$-valued conditional logic **ŁCR** semantically. We denote the language of Łukasiewicz $m$-valued logic as $\mathcal{L}$. $\mathcal{L}$ is given as follows: there is a denumerable set of propositional variables $\Pi$, a unary connective $\neg$, a binary connective $\rightarrow$, and the parentheses. The language of **ŁCR**, $\mathcal{L}_{\succ}$, is obtained by adding the binary connective $\succ$ to $\mathcal{L}$ for conditionals. The set of formulae $\Phi$ is defined recursively:

1. If $p \in \Pi$, $p$ is a formula;
2. If $\phi$ is a formula, so is $\neg\phi$;
3. If $\phi$ and $\psi$ are formulae, so are $(\phi \succ \psi)$ and $(\phi \rightarrow \psi)$.

We agree that the brackets on the outermost side of the formula can be omitted without causing confusion. We use $p$, $q$, … for propositional variables and $\phi$, $\psi$, $\gamma$, $\chi$, … for formulae. Atomic formulae of $\mathcal{L}_{\succ}$ are propositional variables. From the semantic point of view, the set $\{0, \frac{1}{m-1}, \frac{2}{m-1}, \ldots, \frac{m-2}{m-1}, 1\}$, denoted by $\mathcal{T}$, is the truth value set with natural order, and we use $a$, $b$, $c$ and so on to represent the elements in it. Łukasiewicz $m$-valued propositional logic is denoted by **Ł** in this paper.

In $\mathcal{L}$ we can introduce other logical connectives and formulae as follows (cf. [20, 19, 10 (Section 6.2, p.109; Section 9.1.2, p.186; Section 9.1.3, p.194)]):

$$\mathtt{t} =_{\mathrm{df}} p \rightarrow p,$$

$$\mathtt{f} =_{\mathrm{df}} \neg(p \rightarrow p),$$

$$\phi \vee \psi =_{\mathrm{df}} ((\phi \rightarrow \psi) \rightarrow \psi),$$

$$\phi \wedge \psi =_{\mathrm{df}} \neg(\neg\phi \vee \neg\psi),$$

$$\phi \oplus \psi =_{\mathrm{df}} \neg\phi \rightarrow \psi,$$

$$\phi \odot \psi =_{\mathrm{df}} \neg(\phi \rightarrow \neg\psi),$$

$$\phi \ominus \psi =_{\mathrm{df}} \phi \odot \neg\psi,$$

$$\phi \leftrightarrow \psi =_{\mathrm{df}} (\phi \rightarrow \psi) \wedge (\psi \rightarrow \phi),$$

$$\twoheadrightarrow_{i=1}^{0}(\phi_i, \psi) =_{\mathrm{df}} \psi,$$

$$\twoheadrightarrow_{i=1}^{k}(\phi_i, \psi) =_{\mathrm{df}} \phi_k \rightarrow (\twoheadrightarrow_{i=1}^{k-1}(\phi_i, \psi)),$$

$$\twoheadrightarrow_{i=1}^{k}(\phi, \psi) =_{\mathrm{df}} \phi \rightarrow (\twoheadrightarrow_{i=1}^{k-1}(\phi, \psi)),$$

$$\prod_{i=1}^{1}\phi_1 =_{\mathrm{df}} \phi_1,$$

$$\prod_{i=1}^{n+1}\phi_i =_{\mathrm{df}} (\prod_{i=1}^{n}\phi_i) \odot \phi_{n+1},$$

$$\sum_{i=1}^{1}\phi_i =_{\mathrm{df}} \phi_1,$$

$$\sum_{i=1}^{n+1}\phi_i =_{\mathrm{df}} (\sum_{i=1}^{n}\phi_i) \oplus \phi_{n+1},$$

$$\mathsf{J}_1(\phi) =_{\mathrm{df}} \prod_{i=1}^{m-1}\phi \quad (=\prod_{i=1}^{m-1}\phi_i \text{, where } \phi_i = \phi,\ 1 \le i \le m - 1).$$

For each truth value $a$ with $\frac{1}{2} \le a \le 1$ we denote by $n(a)$ the largest integer $k$ such that $k \cdot (1 - a) < 1$. We have $1 - n(a) \cdot (1- a) \le 1 - a$. Let

$\mathsf{J}_a(\phi) =_{\mathrm{df}} \mathsf{J}_1(\neg \prod_{i=1}^{n(a)}\phi \leftrightarrow \phi)$ for $n(a) = \frac{a}{1-a}$ and $a \ge \frac{1}{2}$ (i.e., the case where $k$ satisfies $k \cdot (1 - a) + a =1$),

$\mathsf{J}_a(\phi) =_{\mathrm{df}} \mathsf{J}_{n(a)\,\cdot\,(1-a)}(\neg \prod_{i=1}^{n(a)}\phi)$ for $n(a) > \frac{a}{1-a}$ and $a \ge \frac{1}{2}$ (i.e., the case where $k$ satisfies $k \cdot (1 - a) + a > 1$) and $k \cdot (1 - a) <1$),

$\mathsf{J}_a(\phi) =_{\mathrm{df}} \mathsf{J}_{1-a}(\neg\phi)$ for $a < \frac{1}{2}$,

$\mathsf{I}_a(\phi) =_{\mathrm{df}} \mathsf{J}_a(\phi) \vee \ldots \vee \mathsf{J}_b(\phi) \vee \ldots \vee \mathsf{J}_1(\phi),\ a < b$.

## 2.2. Worlds Semantics

Kripke models have been generalized to many-valued logic and intuitionist conditional logic. In this paper, we will give the language $\mathcal{L}_{\succ}$ a Kripke model with the structure needed to interpret the operator $\succ$. Our presentation closely follows that given for 3-valued modal logic in [22], and refers to a variant of the semantics for the conditional logic given in [3, 4, 32]. First, we give the model of Łukasiewicz $m$-valued conditional logic as follows:

**Definition 2.1.** A *Kripke* **ŁCR** *model* is a structure $\mathfrak{I} = \langle W, \mathcal{A}, \{R_X : X \in \mathcal{A}\}, v\rangle$ such that:

1. $W$ is a nonempty set of worlds;
2. $\mathcal{A} \subseteq (\mathscr{P}(W))^m$ is a set of propositions;
3. $R_X$ is a function from $W \times W$ to $\mathcal{T}$, for every $X = (X_1, \ldots, X_m) \in \mathcal{A}$;
4. $v$ is a function from $\Pi \times W$ to $\mathcal{T}$.

The designated value is 1 ("true"). We denote $v(\phi, x)$ as $v_x(\phi)$, and set $|\phi| = (|\phi|_1, \ldots, |\phi|_m)$ (i.e. the proposition expressed by $\phi$), where $|\phi|_i = \{x \in W\colon v_x(\phi) = \frac{i-1}{m-1}\}$, $1 \le i \le m$. Obviously $|\phi| = |\psi|$ iff $|\phi|_i = |\psi|_i$ for all $i$, $1 \le i \le m$. We also have $|\phi|_i \cap |\phi|_j = \varnothing$ for $i \ne j$, and $\bigcup_{i=1}^{m} |\phi|_i = W$. In the classical conditional logic, a proposition is regarded as the set of possible worlds in which it is true. In the many-valued conditional logic, a proposition is regarded as a $m$-tuple of the sets of possible

worlds, so as to distinguish different many-valued propositions.

We denote ~ and $\rightarrowtail$ the functions defined as follows (for $a$, $b \in \mathcal{T}$):

$$\sim a = 1 - a, \quad a \rightarrowtail b = \min\{1, 1 - a + b\}.$$

Given a Kripke **ŁCR** model $\mathfrak{I}$, $v$ is extended by the following recursive conditions ($\forall x \in W$):

5. $v_x(\neg\phi) = \sim v_x(\phi)$;
6. $v_x(\phi \to \psi) = v_x(\phi) \rightarrowtail v_x(\psi)$;
7. $v_x(\phi \succ \psi) = \bigwedge \{xR_{|\phi|}y \rightarrowtail v_y(\psi): y \in W\}$

In classical conditional logic, $\phi \succ \psi$ is true in world $x$ if and only if $\psi$ is true in all worlds $y$ such that $xR_{|\phi|}y = 1$, i.e. $v_x(\phi \succ \psi) = \bigwedge \{xR_{|\phi|}y \rightarrowtail v_y(\psi): y \in W\}$. We extend this semantics to the many-valued case, and it makes sense: $\phi \succ \psi$ is true in possible worlds $x$ if and only if the true value of $\psi$ is greater than or equal to $a$ in world $y$ such that $xR_{|\phi|}y = a$. This generalizes classical conditional semantics.

By Kripke **ŁCR** model, for arbitrary $v_w(\phi)$, $v_w(\psi) \in \mathcal{T}$, $w \in W$ we have

$v_w(\phi) \wedge v_w(\psi) = \min\{v_w(\phi), v_w(\psi)\}$,

$v_w(\phi) \vee v_w(\psi) = \max\{v_w(\phi), v_w(\psi)\}$,

$v_w(\phi) \oplus v_w(\psi) = \min\{1, v_w(\phi) + v_w(\psi)\}$,

$v_w(\phi) \odot v_w(\psi) = \max\{0, v_w(\phi) + v_w(\psi) - 1\}$,

$v_w(\phi) \ominus v_w(\psi) = \max\{0, v_w(\phi) - v_w(\psi)\}$,

$$v_w(\mathrm{J}_a(\phi)) = \begin{cases} 1, & \text{if } v_w(\phi) = a \\ 0, & \text{if } v_w(\phi) \neq a \end{cases}.$$

For convenience of expression, here the algebraic operations $\oplus$, $\odot$, and $\ominus$ and their corresponding logical connectives are represented by the same symbol. Obviously, $\langle \mathcal{T}, \oplus, \sim, 0 \rangle$ constitutes a linearly ordered $MV_m$-algebra.

Now we give some definitions and construct the Łukasiewicz $m$-valued conditional logic system **ŁCR** which allows the accessibility relation $R_X$ is many-valued.

**Definition 2.2.** Let $\Sigma$ be a set of formulae and $\mathcal{C}$ is a class of Kripke **ŁCR** models, $\Sigma \vDash_{\mathcal{C}} \phi$ if and only if for all worlds $w$ of all $\mathfrak{I} \in \mathcal{C}$, if $v_w(\psi) = 1$ for each $\psi \in \Sigma$, then $v_w(\phi) = 1$. If $\Sigma \vDash_{\mathcal{C}} \phi$, the inference is called *valid* (in $\mathcal{C}$); if, in addition, $\Sigma = \varnothing$, the formula $\phi$ is called *valid* (in $\mathcal{C}$). The class of all Kripke **ŁCR** models will be denoted by $\mathcal{C}_{\mathbf{K}}^{\mathbf{ŁCR}}$.

**Definition 2.3.** A Kripke **ŁCR** model is said to be *finite* just in case $W$ is. A set of formulae (a logic) is *determined by* a class of frames exactly when the set is both sound and complete with respect to the class.

## 2.3. Axiom System **ŁCR**

In this section, we introduce an axiom system for **ŁCR** and establish a few elementary results. **ŁCR** is built on top of an axiomatization of Łukasiewicz $m$-valued propositional logic **Ł**. Rosser and Turquette [20], Tokarz [30], Schotch et al. [21], Grigolia [11], and Tuziak [31] have given the axiomatization of Łukasiewicz $m$-valued logic. **ŁCR** is constituted by **Ł** (all axioms and the rule (MP)) plus the following axioms (axiom schemata) and inference rules:

A1: $(\phi \succ (\psi \wedge \theta)) \rightarrow ((\phi \succ \psi) \wedge (\phi \succ \theta))$ (CM)

A2: $((\phi \succ \psi) \wedge (\phi \succ \theta)) \rightarrow (\phi \succ (\psi \wedge \theta))$ (CC)

A3: $\phi \succ \mathtt{t}$ (CN)

$$\frac{(\phi \leftrightarrow \psi)}{(\phi \succ \theta) \leftrightarrow (\psi \succ \theta)} \quad \text{(RCEA)}$$

$$\frac{(\phi \leftrightarrow \psi)}{(\theta \succ \phi) \leftrightarrow (\theta \succ \psi)} \quad \text{(RCEC)}$$

$$\frac{\twoheadrightarrow_{i=1}^{m} (\mathsf{l}_{a_i \odot b}(\gamma_i), \mathsf{l}_{a \odot b}(\gamma)) \text{ for every } b \in \mathcal{T}}{\twoheadrightarrow_{i=1}^{m} (\mathsf{l}_{a_i}(\phi \succ \gamma_i), \mathsf{l}_{a}(\phi \succ \gamma))}, \text{ where } a \in \mathcal{T}, a_i = \tfrac{m-i}{m-1}, 1 \leq i \leq m \quad (\mathrm{R}_a)$$

Axioms A1-A3 are axioms about conditionals. The reader should notice that **ŁCR** is formed by extending Łukasiewicz $m$-valued propositional logic with the $m$-valued counterparts of the conditional theses in [3]. (CM), (CC), (CN), (RCEA), and (RCEC) correspond respectively to the corresponding axioms or rules in [3]. However, it can be verified that the following conditional K axiom (CK)

$$(\phi \succ (\psi \rightarrow \theta)) \rightarrow ((\phi \succ \psi) \rightarrow (\phi \succ \theta)) \quad \text{(CK)}$$

does not hold in **ŁCR**, so we add the rule ($\mathrm{R}_a$) and call the system **ŁCR** instead of **ŁCK**.

As in standard modal logic, some care must be taken in defining the relation of consequence $\vdash_{\text{ŁCR}}$ that the logic gives rise to. We define it as follows.

**Definition 2.4.** $\vdash_{\text{ŁCR}} \phi$ ($\phi$ is a theorem of **ŁCR**) if and only if $\phi \in$ **ŁCR**. $\Gamma \vdash_{\text{ŁCR}} \phi$ if and only if there is a sequence $\phi_1, \ldots, \phi_n$ such that each formula in $\phi_1, \ldots, \phi_n, \phi$ is a theorem of **ŁCR**, an element of $\Gamma$, or is obtained by (MP) from preceding formulae.

Now let's give some theorems of **Ł**, and they are used directly or indirectly in this paper.

**Proposition 2.5.** Ł *contains the following theorems*:

(1) $\vdash_{Ł}\ \twoheadrightarrow_{i=1}^{m-1}(\phi,\ J_1(\phi))$

(2) $\vdash_{Ł} \phi \wedge \psi \rightarrow \psi \wedge \phi$

(3) $\vdash_{Ł} (J_a(\phi) \rightarrow (J_a(\phi) \rightarrow \psi)) \rightarrow (J_a(\phi) \rightarrow \psi)$

(4) $\vdash_{Ł}\ \twoheadrightarrow_{i=1}^{m}(J_{\frac{i-1}{m-1}}(\phi) \rightarrow \psi, \psi)$

(5) $\vdash_{Ł}\ J_1(\phi) \rightarrow \phi$

(6) $\vdash_{Ł} (\phi \rightarrow (\psi \rightarrow \theta)) \rightarrow (\psi \rightarrow (\phi \rightarrow \theta))$

(7) $\vdash_{Ł} \phi \wedge \psi \rightarrow \phi$

(8) $\vdash_{Ł} \phi \leftrightarrow \neg\neg\phi$

(9) $\vdash_{Ł} \phi \rightarrow (\psi \rightarrow \phi)$

(10) $\vdash_{Ł} (\psi \odot (\psi \rightarrow \theta)) \rightarrow \theta$

(11) $\vdash_{Ł} (\phi \odot \psi \rightarrow \theta) \leftrightarrow (\phi \rightarrow (\psi \rightarrow \theta))$

(12) $\vdash_{Ł} (\phi \rightarrow \psi) \rightarrow ((\phi \rightarrow \theta) \rightarrow (\phi \rightarrow \psi \wedge \theta))$

(13) $\vdash_{Ł} J_a(\psi) \rightarrow \neg I_b(\psi)$, *where* $a < b$.

(14) $\vdash_{Ł}\ \twoheadrightarrow_{i=1}^{m}((J_{\frac{i-1}{m-1}}(\phi) \leftrightarrow J_{\frac{i-1}{m-1}}(\psi)), (\phi \leftrightarrow \psi))$

(15) $\vdash_{Ł} (\phi \leftrightarrow \psi) \rightarrow (J_{\frac{i-1}{m-1}}(\phi) \leftrightarrow J_{\frac{i-1}{m-1}}(\psi))$, $1 \leq i \leq m$.

(16) $\vdash_{Ł}\ J_1(J_a(\psi)) \leftrightarrow J_a(\psi)$

(17) $\vdash_{Ł} (I_a(\phi) \rightarrow (I_a(\phi) \rightarrow \psi)) \rightarrow (I_a(\phi) \rightarrow \psi)$

(18) $\vdash_{Ł}\ I_a(\psi_1 \wedge \psi_2) \leftrightarrow I_a(\psi_1) \wedge I_a(\psi_2)$

(19) $\vdash_{Ł} (I_a(\phi) \rightarrow (I_b(\psi) \rightarrow I_c(\theta))) \leftrightarrow (I_a(\phi) \wedge I_b(\psi) \rightarrow I_c(\theta))$

(20) $\vdash_{Ł}\ J_1(\neg(I_a(\phi))) \leftrightarrow \neg(I_a(\phi))$

(21) $\vdash_{Ł}\ I_a(\psi) \rightarrow I_b(\psi)$, where $a \geq b$.

**Proof:** All these formulas are tautologies, so they are theorems. For the proof of some theorems, see [5 (pp.88-89), 10 (p.109; pp.183-185), 12 (pp.35-46; 65)]. □

**Lemma 2.6.** For every $a_1, \ldots, a_n, a \in \mathcal{T}$,

$$\frac{\twoheadrightarrow_{i=1}^{n} (\mathsf{I}_{a_i \odot b}(\gamma_i),\ \mathsf{I}_{a \odot b}(\gamma)) \text{ for every } b \in \mathcal{T}}{\twoheadrightarrow_{i=1}^{n} (\mathsf{I}_{a_i}(\phi \succ \gamma_i), \mathsf{I}_a(\phi \succ \gamma))} \qquad (\mathrm{R}_a^{a_1 \dots a_n})$$

**Proof:** Suppose $\vdash_{\textbf{ŁCR}} \twoheadrightarrow_{i=1}^{n} (\mathsf{I}_{a_i \odot b}(\gamma_i),\ \mathsf{I}_{a \odot b}(\gamma))$ for every $b \in \mathcal{T}$. When $n < m$, we have $\vdash_{\textbf{ŁCR}} \twoheadrightarrow_{j=1}^{m-n} (\mathsf{I}_{a_j \odot b}(\mathtt{t}),\ \twoheadrightarrow_{i=1}^{n} (\mathsf{I}_{a_i \odot b}(\chi_i),\ \mathsf{I}_{a \odot b}(\gamma)))$ for every $b \in \mathcal{T}$ and all $a_j \in \mathcal{T} \backslash \{a_1, \dots, a_n\}$. By ($\mathrm{R}_a$), $\vdash_{\textbf{ŁCR}} \twoheadrightarrow_{j=1}^{m-n} (\mathsf{I}_{a_j \odot b}(\phi \succ \mathtt{t}),\ \twoheadrightarrow_{i=1}^{m} (\mathsf{I}_{a_i}(\phi \succ \chi_i),\ \mathsf{I}_a(\phi \succ \gamma)))$. Since $\vdash_{\textbf{Ł}} \mathsf{I}_{a_j \odot b}(\phi \succ \mathtt{t})$, we have $\vdash_{\textbf{ŁCR}} \twoheadrightarrow_{i=1}^{n} (\mathsf{I}_{a_i}(\phi \succ \gamma_i),\ \mathsf{I}_a(\phi \succ \gamma))$ by (MP). When $m \le n$, by Proposition 2.5 (18), $\vdash_{\textbf{ŁCR}} \twoheadrightarrow_{i=1}^{k} (\mathsf{I}_{a_i \odot b}(\chi_i),\ \mathsf{I}_{a \odot b}(\gamma))$ for every $b \in \mathcal{T}$, where $k \le m$. According to ($\mathrm{R}_a$) and the first case of $n < m$, $\vdash_{\textbf{ŁCR}} \twoheadrightarrow_{i=1}^{k} (\mathsf{I}_{a_i}(\phi \succ \chi_i),\ \mathsf{I}_a(\phi \succ \gamma))$. By **Ł**, Proposition 2.5 (18) and (19), $\vdash_{\textbf{ŁCR}} \twoheadrightarrow_{i=1}^{n} (\mathsf{I}_{a_i}(\phi \succ \gamma_i),\ \mathsf{I}_a(\phi \succ \gamma))$.

□

A conditional logic is *classical* iff it is closed under the rules of inference (RCEA) and (RCEC); *normal* iff it is closed under the rules (RCEA) and (CK) axiom. Chellas and Weiss give normal classical conditional logic and normal intuitionistic conditional logic [3, 32] respectively; the basic Łukasiewicz $m$-valued conditional logic **ŁCR** we get is *classical* but not normal.

**Definition 2.7.** A set $\Gamma$ of formulae of **Ł** is *syntactically* (*in*)*consistent* if there is no (there is a) formula $\phi$ of **Ł** such that both $\Gamma \vdash_{\textbf{Ł}} \phi$ and $\Gamma \vdash_{\textbf{Ł}} \neg\phi$. And $\Gamma$ is *maximally consistent* if (a) $\Gamma$ is syntactically consistent, and (b) for any formula $\phi$ of **Ł**, $\phi \in \Gamma$ iff $\Gamma \vdash_{\textbf{Ł}} \phi$.

**Proposition 2.8.** *For a set $\Gamma$ of formulae of* **Ł**, *we have* [②]

(1) *If $\Gamma \vdash_{\textbf{Ł}} \phi$ then $\Gamma' \vdash_{\textbf{Ł}} \phi$ for every superset $\Gamma'$* of *$\Gamma$.*

(2) *If $\Gamma \vdash_{\textbf{Ł}} \phi$* then there is a finite subset *$\Gamma'$* of *$\Gamma$ such that $\Gamma' \vdash_{\textbf{Ł}} \phi$.*

(3) *If $\phi$ belong to $\Gamma$, then $\Gamma \vdash_{\textbf{Ł}} \phi$.*

(5) *If $\Gamma \cup \{\phi\} \vdash_{\textbf{Ł}} \psi$ then $\Gamma \vdash_{\textbf{Ł}} \twoheadrightarrow_{i=1}^{m-1} (\phi, \psi)$.*

(6) *If $\Gamma$ is syntactically inconsistent, then $\Gamma \vdash_{\textbf{Ł}} \phi$ for every formula $\phi$.*

(7) *$\Gamma$ is syntactically inconsistent if and only if $\Gamma \vdash_{\textbf{Ł}} \mathtt{f}$.*

(8) *If $\Gamma \cup \{\phi\}$ is syntactically inconsistent, then $\Gamma \vdash_{\textbf{Ł}} \neg \mathsf{J}_1(\phi)$.*

(9) *If $\Gamma \cup \{\neg \mathsf{J}_1(\phi)\}$ is syntactically inconsistent, then $\Gamma \vdash_{\textbf{Ł}} \phi$.*

(10) *If both $\mathsf{J}_a(\phi)$* (a $\in \mathcal{T}$) *and $\mathsf{J}_b(\phi)$ ($b \neq a$) belong to $\Gamma$ then $\Gamma \vdash_{\textbf{Ł}} \mathtt{f}$ and $\Gamma$ is inconsistent.*

② Proposition 2.8 can be applied to **ŁCR**.

**Proof:** The proof is omitted here (cf. [9, 21, 22]). □

**Definition 2.9.** A set $\Gamma$ of formulae of **ŁCR** is *semantically consistent* if there is a truth-value assignment under which all members of $\Gamma$ evaluate to 1 in the Kripke **ŁCR** model $\mathfrak{I} = \langle W, \mathcal{A}, \{R_X : X \in \mathcal{A}\}, v\rangle$.

### 2.4. Soundness

In the following theorem we show that **ŁCR** is sound with respect to $\mathcal{C}_{\mathbf{K}}^{\mathbf{LCR}}$.

**Theorem 2.10. (Soundness)** $\Gamma \vdash_{\mathbf{LCR}} \phi$ *implies* $\Gamma \vDash_{\mathcal{C}_{\mathbf{K}}^{\mathbf{LCR}}} \phi$.

**Proof:** The axioms must be shown to be valid, and the rules must be shown to preserve validity. We examine only those axioms and rules in which $\succ$ occurs.

According to the semantics of conditional connective $\succ$, we take a world $x$ of a Kripke **ŁCR** model $\mathfrak{I}$ to prove that A1-A3 are valid and (RCEA), (RCEC) are validity preserving.

$v_x(\phi \succ (\psi \wedge \theta)) = \bigwedge \{xR_{|\phi|}y \rightarrowtail v_y(\psi \wedge \theta) : y \in W\} = \bigwedge \{xR_{|\phi|}y \rightarrowtail v_y(\psi) : y \in W\} \wedge \bigwedge \{xR_{|\phi|}y \rightarrowtail v_y(\theta) : y \in W\} = v_x((\phi \succ \psi) \wedge (\phi \succ \theta))$. Hence A1 and A2 take the designated value 1.

Take an arbitrary Kripke **ŁCR** model $\mathfrak{I}$, world $x$, $v_x(\phi \succ \mathtt{t}) = \bigwedge \{xR_{|\phi|}y \rightarrowtail v_y(\mathtt{t}) : y \in W\} = 1$, thus A3 is established.

If $\phi \leftrightarrow \psi$, then for an arbitrary world $x$ of a Kripke **ŁCR** model $\mathfrak{I}$, $v_x(\phi) = v_x(\psi)$. Thus, $|\phi| = |\psi|$ and $xR_{|\phi|}y = xR_{|\psi|}y$, then $v_x(\phi \succ \theta) = \bigwedge \{ xR_{|\phi|}y \rightarrowtail v_y(\theta) : y \in W\} = \bigwedge \{ xR_{|\psi|}y \rightarrowtail v_y(\psi \wedge \theta) : y \in W\} = v_x(\psi \succ \theta)$. Therefore, $v_x((\phi \succ \theta) = (\psi \succ \theta))$, and (RCEA) is sound. Similarly, we can easily verify that (RCEC) is sound.

The rule ($\mathrm{R}_a$) is sound. To see this, let $a \in \mathcal{T}$ and $v$ be a valuation from a Kripke model. We assume that $v_x(\mathrm{I}_{a_1 \odot b}(\gamma_1) \rightarrow ( \ldots ( \rightarrow (\mathrm{I}_{a_m \odot b}(\gamma_m) \rightarrow \mathrm{I}_{a \odot b}(\gamma))) \ldots)) = 1$ for every world $x$. We have to prove that $v_x(\mathrm{I}_{a_1}(\phi \succ \gamma_1) \rightarrow ( \ldots ( \rightarrow (\mathrm{I}_{a_m}(\phi \succ \gamma_m) \rightarrow \mathrm{I}_a(\phi \succ \gamma))) \ldots)) = 1$. It suffices to prove that if $v_x(\mathrm{I}_{a_i}(\phi \succ \gamma_i) = 1$ for every $i \in \{1, \ldots, m\}$, then $v_x(\mathrm{I}_a(\phi \succ \gamma)) = 1$. That is, we have to prove that if $a_i \leq v_x(\phi \succ \gamma_i)$ for every $i \in \{1, \ldots, m\}$, then $a \leq v_x(\phi \succ \gamma)$. Hence, let us assume that $a_i \leq v_x(\phi \succ \gamma_i)$ for every $i \in \{1, \ldots, m\}$. We must prove that $a \leq v_x(\phi \succ \gamma)$, i.e. that $a \leq (xR_{|\phi|}y \rightarrowtail v_y(\gamma)$ for every world $y$. Thus, we consider a word $y$ and define $b := xR_{|\phi}y$. If $b \le \neg a$ then it is obvious that $a \le b \rightarrowtail 0 \le (xR_{|\phi}y \rightarrowtail v_y(\gamma)$. Let us now consider the case that $\neg a < b$ (i.e. $0 < a \odot b$). The fact that $a_i \le v_x(\phi \succ \gamma_i)$ tells us that $a_i \odot b \le v_y(\gamma_i)$ for every $i \in \{1, \ldots, \mathrm{m}\}$. Therefore, $v_y(\mathrm{I}_{a_i \odot b}(\gamma_i)) = 1$. Using the assumption about the upper part of the rule ($\mathrm{R}_a$) we get that $v_y(\mathrm{I}_{a \odot b}(\gamma)) = 1$. Thus, $a \le b \rightarrowtail v_y(\phi) = xR_{|\phi}y \rightarrow v_y(\gamma)$ for every world $y$. The proof is finished. □

**2.5. Completeness**

In this section, we prove that **ŁCR** is complete with respect to $\mathcal{C}_{\mathbf{K}}^{\mathbf{ŁCR}}$. The proof proceeds by constructing a canonical model.

In order to build a canonical model, we must first construct the maximum consistent sets of formulae. They are constructed in the usual way [9].

**Proposition 2.11.** *For any formula $\phi$, there is exactly one $a \in \mathcal{T}$ such that, in any maximal consistent set $\Delta$ of formulae, $\mathrm{J}_a(\phi) \in \Delta$. If $\phi \in \Delta$, then it must be that $\mathrm{J}_1(\phi) \in \Delta$.*

**Proof:** If $\forall a \in \mathcal{T}$, $\mathrm{J}_a(\phi) \notin \Delta$, $\Delta \cup \{\mathrm{J}_a(\phi)\}$ is syntactically inconsistent by the maximal consistency of $\Delta$. $\Delta \cup \{\mathrm{J}_a(\phi)\} \vdash_{\mathbf{Ł}} \mathrm{f}$ by Proposition 2.8 (7), and $\Delta \vdash_{\mathbf{Ł}} \rightarrow_{i=1}^{m-1} (\mathrm{J}_a(\phi), \mathrm{f})$ by Proposition 2.8 (5). By repeatedly applying Proposition 2.5 (3) and (MP) rule, we get $\Delta \vdash_{\mathbf{Ł}} \mathrm{J}_a(\phi) \rightarrow \mathrm{f}$. By Proposition 2.5 (4) and (MP) rule, $\Delta \vdash_{\mathbf{Ł}} \mathrm{f}$, a contradiction. Therefore there is at least some $b \in \mathcal{T}$, $\mathrm{J}_b(\phi) \in \Delta$. By Proposition 2.8 (10)，there is exactly an $b \in \mathcal{T}$ such that $\mathrm{J}_b(\phi) \in \Delta$.

Assume $\phi \in \Delta$, we have $\mathrm{J}_1(\phi) \in \Delta$ by Proposition 2.5 (1) and (MP). □

**Proposition 2.12.** *Given a set of formulae $\Gamma$ and a formula $\phi$, if $\Gamma \nvdash_{\mathbf{Ł}} \phi$, then there exists a maximal consistent set $\Delta$ such that $\Gamma \subseteq \Delta$ and $\phi \notin \Delta$.*

**Proof:** Suppose that the Proposition does not hold, i.e., for every maximal consistent set $\Delta$ such that $\Gamma \subseteq \Delta$ and $\Gamma \nvdash_{\mathbf{Ł}} \phi$, $\phi \in \Delta$. By Proposition 2.8 (9), $\Gamma \cup \{\neg \mathrm{J}_1(\phi)\}$ is syntactically consistent. According to the construction of maximal consistent set, there exists $\Delta'$, such that $\neg \mathrm{J}_1(\phi) \in \Delta'$. By Proposition 2.11, there exists a $a \in \mathcal{T}$ such that $\mathrm{J}_a(\phi) \in \Delta'$. If $\phi \in \Delta'$, $\mathrm{J}_1(\phi) \in \Delta'$ by Proposition 2.11, a contradiction. □

**Definition 2.13.** $/\phi/ = (/\phi/_1, \ldots, /\phi/_m)$, $/\phi/_i = \{\Gamma\colon \Gamma$ is maximal consistent set of formulae, and $\mathrm{J}_{\frac{i-1}{m-1}}(\phi) \in \Gamma\}$, $1 \le i \le m$.

**Definition 2.14.** The *canonical model* for **ŁCR** is a structure $\mathfrak{I}^c = \langle W^c, \mathcal{A}^c, \{ R^c_{/\phi/} : /\phi/ \in \mathcal{A}^c\}, V\rangle$ defined as follows:

1. $W^c = \{\Gamma : \Gamma$ is a maximal consistent set of formulae$\}$;
2. $\mathcal{A}^c = \{/\phi/\colon \phi \in \mathcal{L}_{\succ}\}$;
3. For $x, y \in W^c$, $x R^c_{/\phi/} y = \bigwedge\{a \rightarrowtail b\colon\ \mathrm{J}_a(\phi \succ \psi) \in x,\ \mathrm{J}_b(\psi) \in y$, and $a, b \in \mathcal{T}\}$;
4. For $x \in W^c$, $V_x(p) = a$ if and only if $\mathrm{J}_a(p) \in x$, $a \in \mathcal{T}$.

**Proposition 2.15.** *In the caonical model for* **ŁCR** $\mathfrak{I}^c = \langle W^c, \mathcal{A}^c, \{ R^c_{/\phi/} : /\phi/ \in \mathcal{A}^c\}, V\rangle$,

(i) *For $x, y \in W^c$, if $x R^c_{/\phi/} y = b$, then $\{\mathrm{I}_{a\odot b}(\psi)\colon\ \mathrm{J}_a(\phi \succ \psi) \in x$, for every $a \in \mathcal{T}\} \subseteq y$;*

(ii) *For $x, y \in W^c$, if $\{\mathrm{I}_{a\odot b}(\psi)\colon\ \mathrm{J}_a(\phi \succ \psi) \in x$, for every $a \in \mathcal{T}\} \subseteq y$, then $x R^c_{/\phi/} y \ge b$.*

**Proof:** (i) Suppose $x R^c_{/\phi/} y = b$, $\mathrm{J}_a(\phi \succ \psi) \in x$, for some $a \in \mathcal{T}$, but $\mathrm{I}_{a\odot b}(\psi) \notin y$. By Proposition 2.11, for some $c$ such that $c < a \odot b$, $\mathrm{J}_c(\psi) \in y$. By Definition 2.14, $b \le a \rightarrowtail c$, i.e., $a \odot b \le c$, a contradiction.

(ii) Suppose $\{\mathrm{I}_{a\odot b}(\psi)\colon\ \mathrm{J}_a(\phi \succ \psi) \in x$, for every $a \in \mathcal{T}\} \subseteq y$. Then, if $\mathrm{J}_a(\phi \succ \psi) \in x$, for every $a \in \mathcal{T}$, then $\mathrm{J}_c(\psi) \in y$, for some $c \ge a \odot b$. By Definition 2.14, $x R^c_{/\phi/} y = \bigwedge\{a \rightarrowtail c\colon\ \mathrm{J}_a(\phi \succ \psi) \in x,\ \mathrm{J}_c(\psi) \in y$, and $a, c \in \mathcal{T}\} \ge \bigwedge\{a \rightarrowtail a \odot b\colon\ \mathrm{J}_a(\phi \succ \psi) \in x,\ \mathrm{J}_{a\odot b}(\psi) \in y$, and $a, b \in \mathcal{T}\} \ge \mathrm{b}$, where $a \rightarrowtail a \odot b = b \vee {\sim}a$. □

**Lemma 2.16.** *The canonical model for* **ŁCR**, $\mathfrak{I}^c = \langle W^c, \mathcal{A}^c, \{ R^c_{/\phi/} : /\phi/ \in \mathcal{A}^c\}, V\rangle$, *is a well defined Kripke* **ŁCR** *model.*

**Proof:** Two things must be shown: that $R^c_{/\phi/}$ is well defined in Definition 2.14 and that the canonical model satisfies the conditions of Kripke **ŁCR** model.

We prove that the definition is well defined. Suppose $/\phi/ = /\psi/$. Then, $\vdash_{\mathbf{LCR}} \mathrm{J}_a(\phi) \leftrightarrow \mathrm{J}_a(\psi)$ for every $a \in \mathcal{T}$. If not, say $\nvdash_{\mathbf{LCR}} \mathrm{J}_a(\phi) \rightarrow \mathrm{J}_a(\psi)$, then $\mathrm{J}_a(\phi) \nvdash_{\mathbf{LCR}} \mathrm{J}_a(\psi)$ (Because if $\mathrm{J}_a(\phi)$ $\vdash_{\text{ŁCR}}$ $J_a(\psi)$, then $\vdash_{\text{ŁCR}} \twoheadrightarrow_{i=1}^{m-1} (J_a(\phi),\ J_a(\psi))$ by Proposition 2.8 (5). According to Proposition 2.5 (3), $\vdash_{\text{ŁCR}}$ $J_a(\phi) \rightarrow J_a(\psi)$, a contradition). By Proposition 2.12, there is a maximal consistent set $x$ of formulae such that $J_a(\phi) \in x$, but $J_a(\psi) \notin x$, $/\phi/ \neq /\psi/$, a contradiction. By Proposition 2.5 (14), $\vdash_{\text{ŁCR}} \phi \leftrightarrow \psi$, from which it follows by (RCEA) and Proposition 2.5 (15) that $\vdash_{\text{ŁCR}}$ $J_a(\phi \succ \theta) \leftrightarrow J_a(\psi \succ \theta)$ for every $a \in \mathcal{T}$. For $x \in W^c$, by closure, $J_a(\phi \succ \theta) \in x$ if and only if $J_a(\psi \succ \theta) \in x$. By Definition 2.14, thus $R^c_{/\phi/} = R^c_{/\psi/}$.

For the second, note that $W^c$ is nonempty. If $J_a(p) \in x$, $a \in \mathcal{T}$, then $V_x(J_a(p)) = 1$, $V_x(p) = a$. By Proposition 2.11, there is no $b \in \mathcal{T}$ such that $J_b(p) \in x$, $b \neq a$. Therefore, there is at most one $a \in \mathcal{T}$ such that $V_x(p) = a$, $V_x$ is a mapping. □

Below we give the truth lemma of **ŁCR**, the proof of which is rather complicated.

**Lemma 2.17.** **(Truth Lemma)** *If* $\mathfrak{I}^c = \langle W^c, \mathcal{A}^c, \{ R^c_{/\phi/} : /\phi/ \in \mathcal{A}^c\}, V\rangle$ *is the canonical model for* **ŁCR**, *then for all* $\phi$ *and all* $w \in W^c$: $V_x(\phi) = a$ *if and only if* $J_a(\phi) \in x$, $a \in \mathcal{T}$.

**Proof:** In fact, we just need to prove: for all $\phi \in \Phi$ and all $x \in W$, $V_x(\phi) = a$ if $J_a(\phi) \in x$. For the proof of the other direction: since $V_x(\phi) = a$ and $J_a(\phi) \notin x$, by Proposition 2.11, $J_b(\phi) \in x$ for some $b \neq a$. By the first direction we have $V_x(\phi) = b$. But $V_x$ is a function by Lemma 2.16, a contradiction. We only examine the case where formulae are of the form $\phi \succ \psi$; the induction hypothesis is that the lemma is true for $\phi$ and $\psi$. Suppose that $V_x(\phi \succ \psi) \neq a$ ($a \in \mathcal{T}$). It follows that $\exists y(x R^c_{/\phi/} y \rightarrowtail V_y(\psi)) < a)$ or $\forall y(x R^c_{/\phi/} y \rightarrowtail V_y(\psi) \neq a)$, i.e., $\exists y(x R^c_{/\phi/} y \rightarrowtail V_y(\psi)) < a)$ or $\forall y(a < x R^c_{/\phi/} y \rightarrowtail V_y(\psi))$ (obviously $a \neq 1$).

If $\exists y \in W^c$, $x R^c_{/\phi/} y \rightarrowtail V_y(\psi) < a$, i.e., $\exists y \in W^c$, $x R^c_{/\phi/} y \rightarrowtail V_y(\psi) = b$ for some $b$, $b < a$, then $V_y(\psi) = b \odot c$[3], where $x R^c_{/\phi/} y = c$. By the induction hypothesis, $J_{b \odot c}(\psi) \in y$. Since $c \rightarrowtail V_y(\psi) = b$ and $V_y(\psi) \geq 0$, $c \geq \neg b$ and $b \odot c < a \odot c$. By Proposition 2.5 (13), $\neg I_{a \odot c}(\psi) \in y$. According to the maximal consistency of $y$, $I_{a \odot c}(\psi) \notin y$. By Definition 2.14, $J_a(\phi \succ \psi) \notin x$.

If $\forall y \in W^c$, $a < x R^c_{/\phi/} y \rightarrowtail V_y(\psi)$, by the induction hypothesis $I_{(a + \frac{1}{m-1}) \odot b}(\psi) \in y$[4], where $x\, R^c_{/\phi/}\, y = b$. We also need to prove: if $\forall y \in W^c$, $\mathrm{I}_{\left(a+\frac{1}{m-1}\right)\odot b}(\psi) \in y$, then $\mathrm{J}_a(\phi \succ \psi) \notin x$ (i.e. if $\mathrm{J}_a(\phi \succ \psi) \in x$, then $\exists y \in W^c$, $\mathrm{I}_{\left(a+\frac{1}{m-1}\right)\odot b}(\psi) \notin y$). Suppose that $\mathrm{J}_a(\phi \succ \psi) \in x$, then $\exists y \in W^c$ such that $\{\mathrm{I}_{c\odot b}(\gamma)$: $\mathrm{J}_c(\phi \succ \gamma) \in x$, for every $c \in \mathcal{T}\} \subseteq y$, and $\{\mathrm{I}_{c\odot b}(\gamma)$: $\mathrm{J}_c(\phi \succ \gamma) \in x$, for every $c \in \mathcal{T}\}$ $\cup \{\neg \mathrm{I}_{\left(a+\frac{1}{m-1}\right)\odot b}(\psi)\}$ is syntactically consistent by Proposition 2.15 (i), where $x\, R^c_{/\phi/}\, y = b$. If not, i.e. $\forall y \in W^c$, if $\{\mathrm{I}_{c\odot d}(\gamma)$: $\mathrm{J}_c(\phi \succ \gamma) \in x$, for every $c \in \mathcal{T}\} \subseteq y$ (now let $x\, R^c_{/\phi/}\, y = d$), then $\{\mathrm{I}_{c\odot d}(\gamma)$: $\mathrm{J}_c(\phi \succ \gamma) \in x$, for every $c \in \mathcal{T}\} \cup \{\neg \mathrm{I}_{\left(a+\frac{1}{m-1}\right)\odot d}(\psi)\}$ is syntactically inconsistent. Then, by Proposition 2.8 (8), $\{\mathrm{I}_{c\odot d}(\gamma)$: $\mathrm{J}_c(\phi \succ \gamma) \in x$, for every $c \in \mathcal{T}\} \vdash_{\mathbf{ŁCR}} \neg \mathrm{J}_1(\neg \mathrm{I}_{\left(a+\frac{1}{m-1}\right)\odot d}(\psi))$. By Proposition 2.5 (20) and (8), $\{\mathrm{I}_{c\odot d}(\gamma)$: $\mathrm{J}_c(\phi \succ \gamma) \in x$, for every $c \in \mathcal{T}\} \vdash_{\mathbf{ŁCR}} \mathrm{I}_{\left(a+\frac{1}{m-1}\right)\odot d}(\psi)$. By Proposition 2.8 (2), $\{\mathrm{I}_{c_1\odot d}(\gamma_1)$, $\mathrm{I}_{c_2\odot d}(\gamma_2)$, …, $\mathrm{I}_{c_k\odot d}(\gamma_k)\} \vdash_{\mathbf{ŁCR}} \mathrm{I}_{\left(a+\frac{1}{m-1}\right)\odot d}(\psi)$, where $\mathrm{I}_{c_l\odot d}(\gamma_l) \in \{\mathrm{I}_{c\odot d}(\gamma)$: $\mathrm{J}_c(\phi \succ \gamma) \in x$, for every $c \in \mathcal{T}\}$, $1 \le l \le k$. By Proposition 2.8 (5), $\vdash_{\mathbf{ŁCR}} \twoheadrightarrow^{m-1}_{i=1}(\mathrm{I}_{c_1\odot d}(\gamma_1)$, $\twoheadrightarrow^{m-1}_{i=1}(\ldots(\twoheadrightarrow^{m-1}_{i=1}(\mathrm{I}_{c_k\odot d}(\gamma_k), \mathrm{I}_{\left(a+\frac{1}{m-1}\right)\odot d}(\psi)))\ldots))$. By Proposition 2.5 (17) and (6), $\vdash_{\mathbf{ŁCR}}$ $\twoheadrightarrow^{k}_{i=1}(\mathrm{I}_{c_i\odot d}(\gamma_i), \mathrm{I}_{\left(a+\frac{1}{m-1}\right)\odot d}(\psi))$. Since $d$ is arbitrary, we can always get: for every $d_i \in \mathcal{T}$, $1 \le i \le m$,

$$\vdash_{\mathbf{ŁCR}} \twoheadrightarrow^{k_i}_{j=1}(\mathrm{I}_{c^i_j\odot d_i}(\gamma^i_j), \mathrm{I}_{\left(a+\frac{1}{m-1}\right)\odot d_i}(\psi)),$$

where $\mathrm{I}_{c^i_j\odot d_i}(\gamma^i_j) \in \{\mathrm{I}_{c\odot d_i}(\gamma)$: $\mathrm{J}_c(\phi \succ \gamma) \in x$, for every $c \in \mathcal{T}\} \subseteq y_i$. Now suppose we have $\vdash_{\mathbf{ŁCR}} \twoheadrightarrow^{k_1}_{j_1=1}(\mathrm{I}_{c^1_{j_1}\odot d_1}(\gamma^1_{j_1}), \mathrm{I}_{\left(a+\frac{1}{m-1}\right)\odot d_1}(\psi))$ and $\vdash_{\mathbf{ŁCR}} \twoheadrightarrow^{k_2}_{j_2=1}(\mathrm{I}_{c^2_{j_2}\odot d_2}(\gamma^2_{j_2}), \mathrm{I}_{\left(a+\frac{1}{m-1}\right)\odot d_2}(\psi))$, where $\mathrm{I}_{c^1_{j_1}\odot d_1}(\gamma^1_{j_1}) \in y_1$ and $\mathrm{I}_{c^2_{j_2}\odot d_2}(\gamma^2_{j_2}) \in y_2$, according to Proposition 2.5 (9) and Proposition 2.15, we can always get that $\vdash_{\mathbf{ŁCR}} \twoheadrightarrow^{k_2}_{j_2=1}(\mathrm{I}_{c^2_{j_2}\odot d_1}(\gamma^2_{j_2}), (\twoheadrightarrow^{k_1}_{j_1=1}(\mathrm{I}_{c^1_{j_1}\odot d_1}(\gamma^1_{j_1}), \mathrm{I}_{\left(a+\frac{1}{m-1}\right)\odot d_1}(\psi))))$, where $\mathrm{I}_{c^2_{j_2}\odot d_1}(\gamma^2_{j_2}) \in y_1$, and $\vdash_{\mathbf{ŁCR}} \twoheadrightarrow^{k_1}_{j_1=1}(\mathrm{I}_{c^1_{j_1}\odot d_2}(\gamma^1_{j_1}), (\twoheadrightarrow^{k_2}_{j_2=1}(\mathrm{I}_{c^2_{j_2}\odot d_2}(\gamma^2_{j_2}), \mathrm{I}_{\left(a+\frac{1}{m-1}\right)\odot d_2}(\psi))))$, where $\mathrm{I}_{c^1_{j_1}\odot d_2}(\gamma^1_{j_1}) \in y_2$. We can get in this way and by Proposition 2.5 (6): for every $d_i \in \mathcal{T}$, $1 \le i \le m$,

$$\vdash_{\mathbf{ŁCR}} \twoheadrightarrow^{k_1}_{j_1=1}(\mathrm{I}_{c^1_{j_1}\odot d_i}(\gamma^1_{j_1}), \ldots, \twoheadrightarrow^{k_i}_{j_i=1}(\mathrm{I}_{c^i_{j_i}\odot d_i}(\gamma^i_{j_i}), \ldots, \twoheadrightarrow^{k_m}_{j_m=1}(\mathrm{I}_{c^m_{j_m}\odot d_i}(\gamma^m_{j_m}), \mathrm{I}_{\left(a+\frac{1}{m-1}\right)\odot d_i}(\psi))\ldots)\ldots).$$

---

[3] By calculating or referring to the properties of MV-algebra.

[4] The validity of this step can be verified by $a + \frac{1}{m-1} \leq b \rightarrowtail V_y(\psi)$ and the property of MV-algebra ( cf.[5, p.10]).

We can rewrite it as $\vdash_{\mathbf{ŁCR}} \twoheadrightarrow^{\Sigma^m_{l=1}k_l}_{j=1}(\mathrm{I}_{c_j\odot d_i}(\gamma_j), \mathrm{I}_{\left(a+\frac{1}{m-1}\right)\odot d_i}(\psi))$. By Definition 2.14, $\mathrm{I}_{c_j\odot d_i}(\gamma_j) \in \{\mathrm{I}_{c\odot d_i}(\gamma)$: $\mathrm{J}_c(\phi \succ \gamma) \in x$, for every $c \in \mathcal{T}\} \subseteq y_i$, where $x\, R^c_{/\phi/}\, y_i = d_i$. By ($\mathrm{R}^{a_1\cdots a_n}_a$) and A3, $\vdash_{\mathbf{ŁCR}}$ $\twoheadrightarrow^{\Sigma^m_{l=1}k_l}_{j=1}(\mathrm{I}_{c_j}(\phi \succ \gamma_j), \mathrm{I}_{a+\frac{1}{m-1}}(\phi \succ \psi))$. Since $\mathrm{J}_{c_j}(\phi \succ \gamma_j) \in x$, $\mathrm{I}_{c_j}(\phi \succ \gamma_j) \in x$, $1 \le j \le \Sigma^m_{l=1}k_l$. By (MP), $\mathrm{I}_{a+\frac{1}{m-1}}(\phi \succ \psi) \in x$, which contradicts the consistency of $x$ (because $\mathrm{J}_a(\phi \succ \psi) \in x$). Thus $\exists y \in W^c$ such that $\{\mathrm{I}_{c\odot b}(\gamma)$: $\mathrm{J}_c(\phi \succ \gamma) \in x$, for every $c \in \mathcal{T}\} \subseteq y$, and $\{\mathrm{I}_{c\odot b}(\gamma)$: $\mathrm{J}_c(\phi \succ \gamma) \in x$, for every $c \in \mathcal{T}\} \cup \{\neg\, \mathrm{I}_{\left(a+\frac{1}{m-1}\right)\odot b}(\psi)\}$ is consistent and can be extended to a maximal consistent set $y'$ and $\mathrm{I}_{\left(a+\frac{1}{m-1}\right)\odot b}(\psi) \notin y'$. By Proposition 2.15 (ii) and Proposition 2.5 (21), $x\, R^c_{/\phi/}\, y' \geq x\, R^c_{/\phi/}\, y = b$ and $\mathrm{I}_{\left(a+\frac{1}{m-1}\right)\odot b'}(\psi) \notin y'$, where $x\, R^c_{/\phi/}\, y' = b'$. Thus we've proven: if $\mathrm{J}_a(\phi \succ \psi) \in x$, then $\exists y \in W^c$, $x\, R^c_{/\phi/}\, y = b$ and $\mathrm{I}_{\left(a+\frac{1}{m-1}\right)\odot b}(\psi) \notin y$. Hence, if $\forall y \in W^c$, $x\, R^c_{/\phi/}\, y = b$ and $\mathrm{I}_{\left(a+\frac{1}{m-1}\right)\odot b}(\psi) \in y$, then $\mathrm{J}_a(\phi \succ \psi) \notin x$. The proof is done. □

**Theorem 2.18. (Completeness)** *$\Gamma \vDash_{\mathfrak{C}^{\mathbf{ŁCR}}_{\mathbf{K}}} \phi$ implies $\Gamma \vdash_{\mathbf{ŁCR}} \phi$.*

**Proof:** Suppose $\Gamma \nvdash_{\mathbf{ŁCR}} \phi$. By Proposition 2.12, there is a maximal consistent set $x$ of formulae such that $\Gamma \subseteq x$ and $\phi \notin x$. By Proposition 2.11, $\mathrm{J}_a(\phi) \in x$ for some $a \in \mathcal{T}$, $a \neq 1$. Where $\mathfrak{I}^c = \langle W^c, \mathcal{A}^c, \{R^c_{/\phi/} : /\phi/ \in \mathcal{A}^c\}, V\rangle$ is the canonical model for **ŁCR**, $x \in W^c$. By Lemma 2.17, for all $\psi \in \Gamma$, $V_x(\psi) = 1$, but $V_x(\phi) = a$. Therefore, $\Gamma \nvDash_{\mathfrak{C}^{\mathbf{ŁCR}}_{\mathbf{K}}} \phi$. □

## 3. Decidability of ŁCR

In this section, we prove the finite model property for **ŁCR**. Decidability follows as a corollary.

**Definition 3.1.** Let $\Sigma$ be a set of formulae closed under subformulae. For a model $\mathfrak{I} = \langle W, \mathcal{A}, \{R_X : X \in \mathcal{A}\}, v\rangle$, we define *the equivalence relation* $\equiv$ on $W$ by:

$$x \equiv y \text{ iff for every } \phi \in \Sigma,\ v_x(\phi) = a \text{ if and only if } v_y(\phi) = a \text{ for all } a \in \mathcal{T}.$$

By $[x]$ we understand the $\equiv$-equivalence class of worlds in $\mathfrak{I}$ generated by $x$. $[X] = ([X]_1, \ldots, [X]_m) \in [\mathcal{A}]$, where $[X]_i \subseteq [W]$. Now a *filtration of $\mathfrak{I}$ through $\Sigma$* is any model $\mathfrak{I}^* = \langle W^*, \mathcal{A}^*, \{R^*_{X^*} : X^* \in \mathcal{A}^* (\mathcal{P}(W^*))^m\}, v^*\rangle$ in which $W^* = [W]$, $X^* = [X] \in [\mathcal{A}] = \mathcal{A}^*$, $v^*$ is a function from $\Pi \times W^*$ to $\mathcal{T}$, and for every $X^* \in \mathcal{A}^*$, $R^*_{X^*}$ is a function from $W^* \times W^*$ to $\mathcal{T}$ satisfying the conditions:

(a) $v^*_{[y]}(p) = v_y(p)$ whenever $p \in \Sigma$;

(b) $[x]\, R^*_{X^*}\, [y] = \bigvee\{xR_Xy : x \in [x], y \in [y]\}$ (where $\bigvee$ is the supremum of the set).

A frame $\mathfrak{F}^* = \langle W^*, \mathcal{A}^*, \{R^*_{X^*} : X^* \in \mathcal{A}^*,\}\rangle$ is a filtration of a frame $\mathfrak{F} = \langle W, \mathcal{A}, \{R_X : X \in \mathcal{A}\rangle$ through $\Sigma$ just in case $\mathfrak{F}^*$ is the frame of a filtration through $\Sigma$ of some model on $\mathfrak{F}$. Note that filtrations are indeed models (or frames), and that a filtration is finite if $\Sigma$ is.

**Theorem 3.2.** *Let $\Sigma$ be a set of formulae closed under subformulae, and let $\mathfrak{I}^*$ be a filtration of a model $\mathfrak{I}$ through $\Sigma$. Then for every $\phi \in \Sigma$ and every x in $\mathfrak{I}$:*

$$v_x(\phi) = a \text{ iff } v^*_{[x]}(\phi) = a,\ a \in \mathcal{T}.$$

*That is*, $[|\phi|] = |\phi|^*$, *for every $\phi \in \Sigma$.*

**Proof:** The proof is by induction, and the only case of interest is that in which the formula in $\Sigma$ has the form $\phi \succ \psi$. So suppose $v_x(\phi \succ \psi) = a$, i.e., $\bigwedge\{xR_{|\phi|}y \rightarrowtail v_y(\psi): y \in W\} = a$. Then $xR_{|\phi|}y \rightarrowtail v_y(\psi) \geq a$ for all $y \in W$, and $\exists\ z \in W$, $xR_{|\phi|}z \rightarrowtail v_z(\psi) = a$. By the induction hypothesis and clause (b) of Definition 3.1, $[x]\,R^*_{|\phi|^*}\,[y] \rightarrowtail v^*_{[y]}(\psi) \geq a$ and $[x]\,R^*_{|\phi|^*}\,[z] \rightarrowtail v^*_{[z]}(\psi) = a$. Hence, $v^*_{[x]}(\phi \succ \psi) = a$.

For the reverse, assume $v^*_{[x]}(\phi \succ \psi) = a$. Then for all $[y] \in [W]$, $[x]\,R^*_{|\phi|^*}\,[y] \rightarrowtail v^*_{[y]}(\psi) \geq a$ and $\exists\,[z] \in [W]$, $[x]\,R^*_{|\phi|^*}\,[z] \rightarrowtail v^*_{[z]}(\psi) = a$. By the induction hypothesis and clause (b), $\forall\, y \in W$, $xR_{|\phi|}y \rightarrowtail v_y(\psi) \geq a$, and $\exists z \in W$, $xR_{|\phi|}z \rightarrowtail v_z(\psi) = a$. Therefore, $v_x(\phi \succ \psi) = \bigwedge\{\mathrm{x}R_{|\phi|}y \rightarrowtail v_y(\psi)\} = a$.

□

**Proposition 3.3.** *Let $\Sigma$ be a finite set of formulae closed under subformulae. For arbitrary model $\mathfrak{I}$, if $\mathfrak{I}^*$ be a filtration of a model $\mathfrak{I}$ through $\Sigma$, then it contains at most $m^n$ worlds (where n denotes the size of $\Sigma$).*

**Proof:** The worlds of $\mathfrak{I}^*$ are the ≡-equivalence classes in $[W]$. Let $g$ be the function with domain $[W]$ and range $\mathcal{T}^{\Sigma}$ defined by $g([x]) = v^*_{[x]}$, where $v^*_{[x]}$ is the evaluation from $\Sigma$ to $\mathcal{T}$. It follows from the definition of the equivalence relation ≡ that $g$ is well defined and injective. Thus the size of $[W]$ is at most $m^n$, where $n$ is the size of $\Sigma$. □

**Theorem 3.4. (Finite Model Property)** *If $\nvdash_{\mathbf{LCR}} \phi$, there is a finite Kripke* **ŁCR** *model $\mathfrak{I}^{\Gamma} = \langle W^{\Gamma}, \mathcal{A}^{\Gamma}, \{R^{\Gamma}_{X^{\Gamma}} : X^{\Gamma} \in \mathcal{A}^{\Gamma}\}, v^{\Gamma}\rangle$ such that $\nvDash^{\mathfrak{I}^{\Gamma}} \phi$.*

**Proof:** Assume $\nvdash_{\textbf{ŁCR}} \phi$. By Theorem 2.18, $\nvDash_{\mathfrak{C}_K^{ŁCR}} \phi$. Then, there is a canonical model $\mathfrak{I} = \langle W, \mathcal{A}, \{R_X : X \in \mathcal{A}\}, v\rangle$ such that $\nvDash^{\mathfrak{I}} \phi$. Let $\Gamma$ contain $\phi$ and be closed under subformulae; then $\Gamma$ is a finite set. Let $\mathfrak{I}^{\Gamma} = \langle W^{\Gamma}, \mathcal{A}^{\Gamma}, \{ R^{\Gamma}_{X^{\Gamma}} : X^{\Gamma} \in \mathcal{A}^{\Gamma} \}, v^{\Gamma}\rangle$ be a filtration of $\mathfrak{I}$ through $\Gamma$. By Theorem 3.2, $v(\phi) = v^{\Gamma}(\phi)$, it follows that $\nvDash^{\mathfrak{I}^{\Gamma}} \phi$. By Proposition 3.3, $\mathfrak{I}^{\Gamma}$ is finite. □

**Corollary 3.5. ŁCR** *is decidable.*

**Proof:** By Theorem 3.4, **ŁCR** has the finite model property. Moreover, it is axiomatizable using finitely many schemas. Therefore, **ŁCR** is decidable. □

Now we have built a basic conditional logic system **ŁCR** based on Łukasiewicz *m*-valued propositional logic **Ł**. The conditional of this system can be seen as Lewis's fixed strict conditional. Lewis believed that counterfactual conditional should be variable strict conditional. To describe the variable strict conditional, it is necessary to restrict its model and add the corresponding axioms. However, due to the complexity of the many-valued accessibility relation, we cannot directly obtain the variable strict conditional logic system by restricting its model. We can only make other constraints to expand **ŁCR**.

## 4. Extending ŁCR

In the specification of Kripke **ŁCR** model, no significant constraints are imposed on $R_X$. Chellas and Weiss have expanded their conditional logic system respectively [3, 32]. In **ŁCR**, we can add the following constraint to the relationship $R_X$ of the model $(x, y \in W;\ X \in \mathcal{A})$:

$$\text{If } xR_Xy = \tfrac{i-1}{m-1}\text{, then } y \in X_j\text{, for some } j \geq i \qquad \text{(łid)}$$

According to (łid), for any formula $\phi$, if $xR_{|\phi|}y \geq \frac{1}{m-1}$, then $v_y(\phi) = 1$. In classical conditional logic, the choice function $f(\phi, w)$ picks the possible worlds where the antecedent $\phi$ is true, and now we obtain a accessibility relation that meets this requirement. The accessible worlds in classical conditional logic are the possible worlds most similar to the current world, and the accessibility relation takes the value 1. And now $xR_Xy$ takes $m$ truth values, which allow us to consider not only the closest possible world where the antecedent is true, but also the possible world with $a$ ($0 < a \in \mathcal{T}$) degree of similarity to (belief in) the current world. $xR_{|\phi|}y = a$ means that the degree of similarity (or belief) of the possible world $y$ to (in) the world $x$ in which the true value of $\phi$ is greater than or equals to $a$.

The semantic constraint corresponds in the obvious way to the following axiom: for every $a$

$\in \mathcal{T}$,

$$\phi \succ \phi \qquad \text{(ŁID)}$$

The soundness of (ŁID) is obvious. In order to make the extended system complete, we only need to prove its canonicity.

**Proof of canonicity:** Assuming (ŁID) holds, we need to prove: For any maximal consistent set $x$ and the formula $\phi$, if $xR_{|\phi|}y = \frac{i-1}{m-1}$, then $y \in |\phi|_j$, for some $j \geq i$, i.e. $\mathsf{I}_{\frac{i-1}{m-1}}(\phi) \in y$. If $xR_{|\phi|}y = \frac{i-1}{m-1}$ and $\phi \succ \phi \in x$, then $\bigwedge\{a \rightarrowtail b\text{: } \mathsf{J}_a(\phi \succ \psi) \in x,\ \mathsf{J}_b(\psi) \in y, \text{ and } a, b \in \mathcal{T}\} = \frac{i-1}{m-1} \leq \bigwedge\{1 \rightarrowtail b\text{: } \phi \succ \phi \in x,\ \mathsf{J}_b(\phi) \in y\}$ by Definition 2.14. Therefore, when $\mathsf{J}_b(\phi) \in y$, $\frac{i-1}{m-1} \leq b$. Thus we have $\mathsf{I}_{\frac{i-1}{m-1}}(\phi) \in y$, for some $j \geq i$. □

According to Lewis [17], in order to ensure that conditionals are variable strict conditionals, we also need to add two constraints to the semantic model as follows:

If $f(\phi, x) \subseteq |\psi|$ and $f(\psi, x) \subseteq |\phi|$, then $f(\phi, x) = f(\psi, x)$;

Either $f(\phi \vee \psi, x) \subseteq |\phi|$ or $f(\phi \vee \psi, x) \subseteq |\psi|$, then $f(\phi \vee \psi, x) = f(\phi, x) \cup f(\psi, x)$.

These constraints can be translated into:

If $R_X(x) \subseteq Y$ and $R_Y(x) \subseteq X$, then $R_X(x) = R_Y(x)$;

Either $R_{X \cup Y}(x) \subseteq X$ or $R_{X \cup Y}(x) \subseteq Y$, then $R_{X \cup Y}(x) = R_X(x) \cup R_Y(x)$.

However, if we generalize in this way, these conditions are not canonical, so we cannot prove its completeness. However, since we can interpret the value of the relation $xR_Xy$ as the degree of similarity, we do not need to further construct variable strict conditionals. According to the basic conditional logic constructed by Chellas, the system **ŁCR** is basic Łukasiewicz $m$-valued conditional logic, which is also the fixed strict conditional logic understood by Lewis, so we call it strict Łukasiewicz $m$-valued conditional logic.

## 5. Concluding Remarks

In this paper we construct basic Łukasiewicz $m$-valued conditional logic **ŁCR** and its extended system. System **ŁCR** generalizes the accessibility relation to $m$-valued. We have proved the meta-logical properties of **ŁCR** using operators $\mathsf{J}$ and $\mathsf{I}$ by extending propositions from the set of possible worlds to the $m$-tuple of the sets of possible worlds. **ŁCR** has the finite model property.

Although we have established the strict Łukasiewicz $m$-valued conditional logic **ŁCR**, this examination has left a number of interesting subjects unexplored. We only construct a conditional logic for Łukasiewicz $m$-valued logic whose accessibility relation for Kripke models is $m$-valued, and we do not know how to construct a Łukasiewicz conditional logic whose accessibility relation is $n$-valued for arbitrary $n \neq m$. We can also construct other Łukasiewicz $m$-valued conditional logic systems. How to develop the conditional logic of Łukasiewicz infinite-valued propositional logic needs further research, and the corresponding counterparts of Lewis's logic VC and Stalnaker's C2 are also to be investigated. Many-valued conditional logic can be studied not only by means of Kripke semantics, but also by means of other semantics, including the algebraic semantics [14], the topological semantics, etc. All these topics remain for further study.